\documentclass[journal,draftcls,onecolumn,letterpaper,twoside,11pt]{IEEEtran}

\RequirePackage[colorlinks,linkcolor=black,citecolor=black,urlcolor=blue]{hyperref}

\usepackage{color}
\usepackage{verbatim}
\usepackage{amssymb}
\usepackage{cite}

\usepackage{float,latexsym,amsfonts,amstext,times}

\ifCLASSINFOpdf
  \usepackage[pdftex]{graphicx}
  \graphicspath{{./Images/}}
\else
\fi

\usepackage[export]{adjustbox}

\usepackage{amsmath}
\usepackage{amsthm}
\usepackage{mathrsfs} 
\usepackage{dsfont}  
\theoremstyle{definition}

\newtheorem{theorem}{Theorem}
\newtheorem{corollary}{Corollary}[theorem]
\DeclareMathOperator*{\argmax}{arg\,max}

\newtheorem{lemma}[theorem]{Lemma}

\newcommand{\Pro}{{\mathbb{P}}}
\newcommand{\Exp}{{\mathbb{E}}}
\usepackage{array}
 \usepackage[caption=false,font=normalsize,labelfont=sf,textfont=sf]{subfig}
\usepackage{fixltx2e}

\DeclareMathOperator*{\esssup}{esssup}

\usepackage{color, colortbl}
\definecolor{Gray}{gray}{0.9}
\newcolumntype{g}{>{\columncolor{Gray}}c}

\usepackage{algorithm,algorithmic}
\newcommand{\bR}{\mathbb{R}}
\newcommand{\bN}{\mathbb{N}}
\newcommand{\bS}{\mathbb{S}}

\newcommand{\cJ}{\mathcal{J}}
\newcommand{\cI}{\mathcal{I}}
\newcommand{\cC}{\mathcal{C}}
\newcommand{\cN}{\mathcal{N}}

\newcommand{\cF}{\mathcal{F}}

\begin{document}

\title{Worst-Case Misidentification Control in Sequential Change Diagnosis using the min-CuSum}

\author{Austin~Warner
and        Georgios~Fellouris
\thanks{This research was supported by the US National Science Foundation under grant AMPS 1736454 through the University of Illinois Urbana-Champaign. Some parts of this paper were presented at the International Symposium on Information Theory (ISIT) 2022 \cite{warner2022cusum}.}
\thanks{A. Warner and G. Fellouris are with the Department of Statistics of the University of Illinois Urbana-Champaign. Email: {awarner5,fellouri}@illinois.edu.}}


\maketitle

\begin{abstract}
The problem of sequential change diagnosis is considered, where a sequence of independent random elements is accessed sequentially, there is an abrupt change in its distribution at some unknown time, and there are two main operational goals: to quickly detect the change and, upon stopping, to accurately identify the post-change distribution among a finite set of alternatives. The  algorithm that raises an alarm as soon as the CuSum statistic that corresponds to one of the post-change alternatives exceeds a certain threshold is studied. When 
the data are generated over independent channels and the change  can occur in only one of them, its worst-case with respect to the change point conditional probability of misidentification,  given that there was not a false alarm,  is shown to decay exponentially fast in the threshold. As a corollary,
in this setup, this algorithm is shown to asymptotically minimize  Lorden’s detection delay criterion, simultaneously  for every possible post-change distribution, within the class of schemes that satisfy prescribed bounds on the   false alarm rate and the worst-case conditional probability of misidentification,  as the former goes to zero sufficiently faster than the latter.  Finally, these theoretical results are also illustrated in  simulation studies.
\end{abstract}




%

\section{Introduction}
%
%
%
%

\IEEEPARstart{T}{he} problem of quickly  detecting a change in  sequentially acquired data dates back to the  pioneering works of Shewhart~\cite{shewhart1931economic} and Page~\cite{page1954continuous}  and is motivated by a wide range of engineering and scientific applications. Examples of 
 such applications can be found in  industrial process quality control~\cite{bissell1969cusum,hawkins2003changepoint,joe2003statistical}, target detection and identification~\cite{ru2009detection,blackman2004multiple},  integrity monitoring of navigation systems~\cite{nikiforov1993application,bakhache1999reliable},
 target tracking~\cite{tartakovsky2003sequential }, network intrusion detection~\cite{tartakovsky2006detection,tartakovsky2006novel}, bioterrorism~\cite{rolka2007issues,fienberg2005statistical},  genomics~\cite{siegmund2013change}.   In most of these applications,  there are typically many  possible types of change and  it is useful, if not critical, to not only detect the change quickly, but also to correctly identify  the type of change upon stopping.  The problem of simultaneously detecting a change in the distribution of  sequentially collected data and  identifying  the correct post-change distribution among a finite set of   alternatives is  known as \textit{sequential change diagnosis}.

One of the first algorithms for the problem of sequential change detection, where the identification aspect is not present, was  Page's Cumulative Sum (CuSum) algorithm~\cite{page1954continuous}. In the case of iid data before and after the change  with completely specified  distributions,   this algorithm admits a recursive structure, which is very  convenient for its practical implementation. Moreover, it minimizes  \textit{Lorden's  criterion}\cite{lorden1971procedures}, i.e., the worst-case conditional expected detection delay with respect to both the change point and the observations up to it,  subject to a user-specified bound on the false alarm rate \cite{moustakides1986optimal}. 

When the post-change distribution is not completely specified but it is assumed to belong to a finite collection of plausible distributions, one can  compute  the  CuSum statistic for each possible post-change distribution and declare that the  change  has occurred as soon as one of them becomes larger than some threshold~\cite{page1954continuous, barnard1959control}. This algorithm, referred to in this work as the min-CuSum, achieves the optimal detection delay, in Lorden's sense,  to (at least) a first-order asymptotic approximation as the false alarm rate goes to zero~\cite{lorden1971procedures}.  It has been studied, in particular, in the case of the \textit{two-sided problem}, where  there is either an increase or decrease in the mean of the observations  after the change (e.g.,~\cite{khan1984cumulative,dragalin1997design, hadjiliadis2006optimal}), as well as in the case of the   \textit{multichannel  problem}, where there are several independent  data sources, e.g., a sensor network, and the change may occur in one or more of them \cite{fellouris2016second,tartakovsky2005asymptotic,tartakovsky2006detection, tartakovsky2006novel,tartakovsky2008asymptotically,hadjiliadis2009one}.




 The first scheme that was proposed for the sequential change \textit{diagnosis} problem
was based on a  generalization of the CuSum algorithm, and termed the ``Generalized CuSum'' \cite{nikiforov1995generalized}. However,  unlike the original CuSum,  this algorithm does not  admit a recursive structure even in the iid setting and its required number of operations  per time instant grows linearly with the number of \textit{all} already collected observations.  Since the change can take a very long time to occur,  the Generalized CuSum  needs to be modified to be applicable in practice. One such modification is to restrict the required number of operations at each time instant, as  suggested in~\cite{lai2000sequential}. A different modification  was  introduced in~\cite{oskiper2002online},  according to which an alarm is raised   as soon as the  CuSum statistics of a  post-change hypothesis against the pre-change and all other post-change  alternatives exceed a certain threshold. The latter algorithm  is referred to as the ``Matrix CuSum,'' as it  requires the  parallel computation of a matrix of CuSum statistics.   A modification of the latter
was proposed in~\cite{warner2022sequential}, which resets  the Matrix CuSum statistics according to an adaptive window. A recursive algorithm which requires only the computation of the original min-CuSum statistics 
 was introduced in~\cite{nikiforov2000simple}. According to it,  an alarm is raised as soon as the CuSum statistic  of a certain post-change alternative against the pre-change exceeds the corresponding   statistic of any  other post-change alternative by a certain threshold. Moreover, the sequential change diagnosis problem has been considered in the special case of the  multichannel problem, but under general stochastic models in \cite{tartakovsky2021asymptotic,pergamenchtchikov2022minimax}.
 



Unlike the  min-CuSum, which requires a number of operations that is linear  in the number of post-change hypotheses, the 
above algorithms  require a number of calculations per time instant that is quadratic in the number of post-change hypotheses.  This can be a disadvantage in certain applications where the change diagnosis problem is relevant, as  for example  in the case of a relative large sensor network.
It is also worth noting that, in the latter context, the min-CuSum admits a decentralized implementation~\cite{hadjiliadis2009one}, 
i.e., it requires the sensors to communicate to a fusion center only to report that they have  raised an alarm.  On the contrary, all the above change diagnosis algorithms  require that,  at each time instant, each sensor communicate its data continuously to the fusion center.

The above practical considerations, together with the  strong theoretical support in the context  of the  sequential \textit{change detection} problem,  can  explain why  the min-CuSum is often used in practice as an ad-hoc solution to the sequential \textit{change diagnosis} problem~\cite{chen2015quickest,wang2021detection, yang2020control}. However, the theoretical understanding of its performance as a change diagnosis scheme  is relatively limited. Its identification rule was shown to be   a strongly consistent estimator of the true post-change distribution in~\cite{han2007detection}, whereas  its probability of misidentification  was studied  in the two-sided problem in 
 ~\cite{khan1981note}, and was claimed to   decay exponentially as a function of the tuning threshold  in a special case of the multichannel problem  in  \cite[Chapter 9.2]{tartakovsky2014sequential}. However, all these results and  statements refer  to the case where  \textit{the change occurs as soon as monitoring begins}.  While this corresponds to the worst-case scenario for the detection delay, this is not the case for  the task of  identification. Thus, the above results do not  provide any guarantees regarding the  conditional probability of wrong identification given that there has not been a false alarm by the time of the change.

Our first contribution in the present paper is that, in the context of the   general sequential change diagnosis problem, we bound the conditional probability of  false identification, given that there has not been a false alarm by  the time of the change. Specifically, we show, under certain conditions, 
that  this probability decays at most exponentially / quadratically in the threshold
if the post-change distributions are ``closer'' / ``as close'' to the pre-change distribution than to each other (in the sense of Kullback-Leibler divergence).   

Based on this general result we show that, in the multichannel problem in which the change can occur only  in one channel, 
 the \textit{worst-case} (with respect to the  change-point) conditional probability of false identification 
decays exponentially fast in the threshold. As a result, we show that  the min-CuSum  asymptotically minimizes  Lorden’s criterion, simultaneously  for every possible post-change distribution, within the class of schemes that satisfy prescribed bounds on the   false alarm rate and the  \textit{worst-case} (with respect to the  change-point)  conditional probability of misidentification,  as the former goes to zero sufficiently faster than the former.

The remainder of the paper is organized as follows. In Section~\ref{sec:state}, we formulate  the problem of sequential change diagnosis.
In Section~\ref{sec:CUSUM},
we describe the min-CuSum algorithm,  review its properties, and  provide intuition on its performance in the sequential change diagnosis problem.  In Section~\ref{sec:results},
we state the theoretical results of the paper.  In Section~\ref{sec:simulation}, we examine the conditional misidentification probability  of the min-CuSum in two simulation studies, and in Section~\ref{sec:conclusion}  we conclude.
 Appendixes~\mbox{\ref{app:cusums}--\ref{app:lemmas}} contain proofs of the main and supporting results.


We end this section by introducing some notation that we use throughout the paper. Thus,
$\bN$ and $\bN_0$ are the sets of 
of positive and non-negative integers, respectively,  and $\bR$ is the set of real numbers. 
For~$x \in \bR$,~$x^+$ is its positive part, i.e.,~$x^+ \equiv \max \lbrace x, 0 \rbrace.$ 
We denote by $x \wedge y$ the minimum of $x, y$. If $(x_n), (y_n)$ are two sequences of positive numbers, then $x_n \sim y_n$ stands for $x_n/y_n \to 1$ as $n \to \infty$. 
The indicator function is written as $\mathds{1}(\cdot)$.

\section{Problem Statement} \label{sec:state}
Let  $X\equiv  \{X_n, \, n \in \bN\}$  be a  sequence of  independent \mbox{$\bS$-valued}  random elements, where $\bS$ is an arbitrary Polish space. For each $n \in \bN$,  we denote by $\cF_n$ the $\sigma$-algebra generated by the first $n$ observations, i.e.,  
$\cF_n \equiv  \sigma(X_1, \ldots, X_n)$, and by $\cF_0$  the trivial $\sigma$-algebra. 

 We assume that   each term of $X$  has  a  density  with respect to a $\sigma$-finite measure $\lambda$, which  is  $f$ up to and including  some deterministic time $\nu \in \bN_0$ and $g$ after $\nu$.  The change point, $\nu$, is completely unknown, but we assume that the   pre-change density, $f$, is completely specified  and that  there is a finite collection of possible post-change distributions, i.e., 
 \begin{equation}
 g \in \{ g_i, \, i \in \mathcal{I} \},
 \end{equation}
 where $\cI$ is an arbitrary finite set.  Moreover, we assume that, for each $i,j \in \cI$ with $i \neq j$, the Kullback-Leibler divergences
  \begin{align} \label{KL}
  \begin{split}
 I_{i} &\equiv  \int_{\bS} \log(g_i/f) \, g_i \, d \lambda, \\ 
  I_{ij} &\equiv \int_{\bS}  \log(g_i/g_j) \, g_i \, d \lambda,  
 \end{split}
\end{align}
are positive and finite. 
 
 We denote by  $ \mathbb{P}_{\infty}$ the  distribution of $X$, and by $\mathbb{E}_{\infty}$ the corresponding expectation, when the change never occurs, i.e., when 
$X$ is a sequence of independent  random elements with common density $f$. Moreover, we denote by  $ \mathbb{P}_{\nu,i} $ the distribution of~$X$, and by $ \mathbb{E}_{\nu,i} $ the corresponding expectation,  when the change occurs at time~$ \nu $ and the post-change density is~$g_i$. For simplicity, when the change occurs from the very beginning,  we suppress the dependence on the change point  and  set~$ \mathbb{P}_i \equiv\mathbb{P}_{0,i}$ and~$\mathbb{E}_i \equiv \mathbb{E}_{0,i}$.

We assume that the terms of the sequence~$X$ are observed sequentially and that the problem  is to  quickly detect the change and, at the same time, to  identify the correct post-change distribution. Thus, we need to specify a random time,~$T$, at which we declare that the change has occurred, and an $\mathcal{I}$-valued random variable,~$D$, that represents our decision regarding the post-change density at the time of stopping. That is, for any $n \in \bN$ and $i \in \mathcal{I} $,  on the event  $\{T=n, D=i \}$  the alarm is raised and   $g_i$  is  declared as the correct post-change density after having taken $n$ observations.  We refer to  such a pair $(T,D)$ as a \textit{diagnosis procedure} if 
 $T$ is  an $\{\cF_n\}$-stopping time and   $D$  is  $\mathcal{F}_T$-measurable, i.e., if 
 $$\{T=n \}, \{T=n, D=i \} \in \cF_n \quad \forall \; n \in \bN, \, i \in \mathcal{I}.$$  
 
Without loss of generality,  we  restrict ourselves to stopping times that are not almost surely bounded under~$\Pro_{\nu,i}$, for any ~$\nu \in \bN$,~$i \in \mathcal{I}$,  and we denote by~$\cC$ the family of all  diagnosis procedures.  

We measure the ability of~$(T,D) \in \cC$   to avoid  false alarms using the average number of  observations  until stopping  under~$\Pro_\infty$, i.e., $\Exp_\infty[T]$, 
and we denote by~$\cC(\alpha)$ the subfamily of diagnosis procedures that control the  expected time until stopping  under~$\Pro_\infty$  by  at least~$1/\alpha $,  i.e., 
 $$\cC(\alpha) \equiv \{(T,D) \in \cC: \Exp_{\infty}[T] \geq 1/\alpha\},$$ 
where $\alpha \in (0,1)$ is a  user-specified value that represents tolerance to false alarms. 
 
To measure the  ability of~$(T,D) \in \cC$  to  identify the  post-change density~$g_i$ when the change occurs at time~$\nu$, we use the conditional probability of  an incorrect identification  given that there was no false alarm, i.e., ~$ \Pro_{\nu, i}(D \neq i \, | \, T > \nu)$, and we denote  by~$\cC(\alpha, \beta)$ the subfamily of~diagnosis procedures in $\cC(\alpha)$ that control the worst-case  conditional probability  of  false identification  below~$ \beta$, i.e., 
\begin{align}\label{def:diagnosis_family}
\left\{ (T,D) \in \cC(\alpha)  :  
\max_{i \in \mathcal{I}} 
\sup_{\nu  \in \bN_0}  \Pro_{\nu,i}(D \neq i \, | \, T > \nu)\leq \beta\right\},
\end{align}
where~$\beta \in (0,1)$ is a  user-specified value that represents tolerance to  false identifications.  

Finally, to measure the  ability of~$(T,D) \in \cC$ to quickly detect the change when the post-change density is~$g_i$  for some~${i \in \mathcal{I}}$,  we adopt  \textit{Lorden's  criterion}~\cite{lorden1971procedures} and employ the worst-case conditional expected detection delay  with respect to both  the change point and  the data up to the change:
  \begin{equation}\label{eqn:lorden_delay_def}
 \cJ_i[T]\equiv \sup_{\nu \in \bN_0} \esssup \Exp_{\nu,i}[ T - \nu \, | \,  \mathcal{F}_{\nu}, T>\nu].
\end{equation}
Our goal in this work is to show that, under certain conditions, the min-CuSum achieves
\begin{equation} \label{infimum}
 \inf_{ (T,D) \in \cC(\alpha, \beta)}  \cJ_i[T],
 \end{equation}
  simultaneously for every~${i \in \mathcal{I}}$, 
to a first-order asymptotic approximation as~$\alpha$ and~$\beta $ go to zero, as long as~$\alpha$ goes to zero sufficiently faster than~$\beta$.

We next provide two setups in which the above sequential change diagnosis problem arises naturally.

\subsection{The Multichannel Setup}\label{subsec:multichannel}
Suppose that $d$ independent data channels are simultaneously monitored and there is a change in the marginal distributions of  an unknown subset of them at some  unknown time,~$\nu$. Specifically,  suppose that channel~$j \in \lbrace 1, \ldots, d \rbrace$ generates a sequence of independent~$\bS_j$-valued random elements, where~$\bS_j$ is some  Polish space, and  let~$X_{j,n}$ denote the observation from channel~$j$ at time~$n$.  If the change does not occur in channel~$j$, then~$X_{j,n}$ has density~$p_j$ with respect to a~$\sigma$-finite measure~$\lambda_j$ on~$\bS_j$ for every~$n \in \bN$\@. On the the other hand,  if the change does occur in that channel at time~$\nu$, then the density of~$X_{j,n}$  is~$p_j$ for~$n \leq \nu$ and $q_j$  for~$n > \nu$. We assume that, for each~$j \in \mathcal{I}$, the densities $p_j$ and~$q_j$ have positive and finite Kullback-Leibler divergences, i.e., 
 \begin{align} \label{multichannel KL conditions}
 \begin{split}
& \int_{\bS_j}\log (q_j/ p_j)  \, q_j \; d\lambda_j \in (0,\infty), \\
& \int_{\bS_j} \log (p_j/q_j)  \, p_j \; d\lambda_j \in (0,\infty).
\end{split}
\end{align} 
This is a special case of the general setup we consider, where 
 \begin{align} \label{multichannel_X}
 X_n &= (X_{1,n}, X_{2,n}, \ldots, X_{d,n}) \in \bS  \equiv \bS_1  \otimes \cdots \otimes \bS_d, 
 \end{align} 
 and the pre-change density is 
 \begin{equation}\label{multichannel_prechange}
  f(x_1, \ldots , x_d) = \prod_{j=1}^d p_j(x_j), \quad   (x_1, \ldots, x_d) \in \bS. 
\end{equation} 
We next consider two special cases for the post-change regime.  In the first one, the change is assumed to   occur in only one channel.  In the second, there is no  prior  information about the number of affected channels, which means that the change may occur in any possible subset of channels. 

\subsubsection{Single-Fault}\label{subsubsec:single_fault}

If the  change can occur in only one channel, as it is often  assumed in the literature (e.g., \cite{pergamenchtchikov2022minimax, tartakovsky2021asymptotic}), then 
$ \mathcal{I} = \{ 1 , \ldots, d \}$ and, for each~$i \in \cI$, ~$g_i$ is  the post-change density when the fault occurs in sensor~$i$, i.e., for $(x_1, \ldots, x_d) \in \bS$\@, \begin{equation} \label{multichannel_single}
 g_i (x_1, \ldots, x_d) = q_i(x_i) \prod_{\substack{j \in \{1 , \ldots, d \} \\ j \neq i}} p_j(x_j).
 \end{equation} 
In this case, the  Kullback-Leibler divergences in~\eqref{KL} take the form
\begin{align} \label{multichannel KL}
\begin{split}
I_i &= \int_{\bS_i} \log (q_i/ p_i)  \, q_i \; d\lambda_i,\\
 I_{ij} &= \int_{\bS_i} \log (q_i/ p_i)  \, q_i \; d\lambda_i + \int_{\bS_j} \log (p_j/ q_j)  \, p_j \; d\lambda_j.
 \end{split}
\end{align}

 \subsubsection{Concurrent-Fault}\label{subsubsec:multiple_fault}

When there is no prior information about the number of  channels that are affected by the change, as for example in~\cite{mei2010efficient, xie2013sequential, fellouris2016second},  then~$\mathcal{I}$ is the family of all non-empty subsets of~$\lbrace 1, \ldots, d \rbrace$ 
and there are~$ 2^d - 1$ possibilities for the post-change density.  Specifically, if   ${i = \lbrace c_1, \ldots, c_p \rbrace}$ where  $c_1, \ldots, c_p \in \lbrace 1, \ldots, d \rbrace$, then $g_i$ represents the post-change density when the  the change occurs in  channels $c_1, \ldots, c_p$, i.e., 
\begin{equation}\label{multichannel_multiple}
g_{i}(x_1 \ldots, x_d) =  \prod_{t \in i } q_t(x_t) \; 
\times \; \prod_{t \notin i } p_t(x_t) ,
\end{equation}
where~$(x_1, \ldots, x_d) \in \bS$, in which the  Kullback-Leibler divergences in~\eqref{KL} take the form
\begin{align}\label{eqn:KL_concurrent}
\begin{split}
I_i &=\sum_{t \in i} \int_{\bS_t} \log(q_t / p_t) q_t \; d \lambda_t  , \\
I_{ij} &=  \sum_{t \in i \setminus j} \int_{\bS_t} \log(q_t / q_t) q_t \; d \lambda_t \,  \\
&+   \sum_{t \in j \setminus i} \int_{\bS_t} \log(p_t / q_t) p_t \; d \lambda_t.
\end{split}
\end{align}


\subsection{The Two-Sided Setup}\label{subsec:two-sided}
Suppose that there is a single stream of observations, whose distibution is specified up to an unknown parameter, and, at some  unknown time, there is either an increase or a decrease in the value of this parameter.  For simplicity,  consider the case that the underlying distribution belongs to a one-parameter exponential family. To be specific, 
let $h$ be a density with respect to $\nu$ and let $\varphi$ denote its cumulant generating function, i.e., 
$$\varphi(\gamma) \equiv \log \int_{\bS} \exp \lbrace \gamma x \rbrace h(x) \, \nu(dx), \quad \gamma \in \bR.$$
Suppose  that the essential domain of $\psi$, $\Gamma$, is non-empty,  i.e.,
\begin{equation}
\Gamma \equiv \{ \gamma \in \bR: \varphi(\gamma) < \infty \} \neq  \emptyset,
\end{equation}
and, for any $\gamma \in \Gamma$, consider the following density
 \begin{equation}\label{def:natural_exp_fam}
h_{\gamma}(x) \equiv h(x) \exp \lbrace \gamma x - \varphi(\gamma) \rbrace,  \quad x \in \bR. 
\end{equation}
Let  $\gamma_0, \gamma_1,  \gamma_2 \in \Gamma$. If the initial parameter is $\gamma_0$ and changes to either $\gamma_1$ or $\gamma_2$, then 
\begin{equation}\label{def:two_sided_densities}
f = h_{\gamma_0}, \quad g_1 = h_{\gamma_1}, \quad g_2 = h_{\gamma_2},
\end{equation}
and the Kullback-Leibler divergences in \eqref{KL}, for $i, j \in \lbrace 1, 2 \rbrace$, take the following form:
\begin{align}
\begin{split}\label{eqn:two-sided_KL1}
 I_i &= \Exp_i\left[\log  \left( \frac{h_{\gamma_i}(X_1)}{h_{\gamma_0}(X_1)}\right) \right] \\
&= (\gamma_i - \gamma_0)  \, \varphi'(\gamma_i) - (\varphi(\gamma_i) - \varphi(\gamma_0))  , \\
\end{split}
\end{align}

\begin{align}
\begin{split}\label{eqn:two-sided_KL2}
I_{ij}&= \Exp_i\left[\log \left(  \frac{h_{\gamma_i}(X_1)}{h_{\gamma_j}(X_1)}\right) \right] \\
&= (\gamma_i - \gamma_j)  \, \varphi'(\gamma_i) - (\varphi(\gamma_i) - \varphi(\gamma_j)).
\end{split} 
\end{align}

\section{The min-CuSum Algorithm}\label{sec:CUSUM}

In this section, we introduce the min-CuSum  algorithm  and provide some intuition regarding its performance in the  sequential change diagnosis problem.  In what follows,  for every  $n \in \bN$ and  $i,j  \in \mathcal{I}$ with $j \neq i$ we set 
\begin{align}\label{LLR}
\begin{split}
 \Lambda_{i}(n) \equiv 
\frac{ g_i (X_n)}{ f(X_n) }, \quad 
  &\ell_{i}(n) \equiv \log  \Lambda_{i}(n), \\
   \Lambda_{ij}(n) \equiv 
\frac{ g_i (X_n)}{ g_j(X_n)}, \quad 
  &\ell_{ij}(n) \equiv \log \Lambda_{ij}(n).
 \end{split}
\end{align}

\subsection{Min-CuSum for Sequential Change Detection}
 Page's CuSum algorithm~\cite{page1954continuous} for detecting a  change from $f$ to $g_i$, for some $i \in \mathcal{I}$,  raises an alarm as soon as   the statistic
\begin{align} \label{cusum_statistic}
Y_{i} (n) &\equiv \max_{ 0 \leq m \leq n}  \sum_{u=m+1}^n  \ell_{i}(u),
\quad n \in \bN,
\end{align} 
exceeds a threshold $b_i>0$,  i.e., at 
\begin{align} \label{cusum}
\begin{split}
\sigma_i (b_i) \equiv \inf \lbrace n \geq 1 : Y_i(n) \geq b_i \rbrace ,
 \end{split}
\end{align}
where we adopt the convention $ \sum_{n+1}^{n} = 0$\@.
An important property of this stopping rule concerning its implementation in practice is that its statistic can be computed via the following recursion:
 \begin{align}\label{cusum_recursion}
Y_{i} (n)  &= \left(Y_{i}(n-1) + \ell_{i}(n) \right)^+, \quad n \in \bN,
\end{align} 
with~$Y_{i}(0) = 0$\@.  Moreover, it  was shown in~\cite{moustakides1986optimal} to   minimize~$\cJ_i$ in~$\cC (\alpha)$, i.e., 
 \begin{align} \label{CUSUM_O}
  \cJ_{i} \left[ \sigma_i(b) \right] =  \inf_{(T, D)  \in \mathcal{C}(\alpha)} \mathcal{J}_i[T]  ,
 \end{align} 
as long as the  false alarm constraint is satisfied with equality.

Of course, this algorithm is directly applicable in our setup  only when $|\mathcal{I}|=1$\@. When~$|\mathcal{I}|>1$, a straightforward way to detect the change  is to   run in parallel the $|\mathcal{I}|$ CuSum statistics, $Y_1, \ldots, Y_{|\mathcal{I}|}$, and  stop as soon as one of them exceeds its corresponding threshold, i.e.,   at $$ \min_{i \in \mathcal{I}} \sigma_i(b_i).$$
In what follows, for the sake of simplicity, we consider equal thresholds, i.e., we set $ b_1 = \cdots = b_{|\mathcal{I}|} = b$, and we refer to the  stopping time
\begin{equation}\label{def:minCuSum}
     \sigma(b) \equiv \min_{i \in \mathcal{I}} \sigma_i(b)
\end{equation}
as  the ``min-CuSum" stopping time. For this stopping time,  it is known (see the proof of~\cite[Theorem 1]{tartakovsky2002efficient}) that, for any~$b>0$,
\begin{align} \label{eqn:CUSUM_ARL}
 \Exp_{\infty}[\sigma(b)] &\geq e^{b}/|\mathcal{I}|,
\end{align} 
 thus,  it is clear that for any~$\alpha \in (0,1)$ we have $\sigma(b_\alpha) \in \cC (\alpha)$ where
\begin{equation}\label{eqn:b_alpha}
 b_\alpha \equiv |\log \alpha| + \log |\mathcal{I}| .
\end{equation}

Moreover, it is well known (see, e.g., ~\cite{lorden1971procedures}) that the worst-case scenario for the conditional expected delay~\eqref{eqn:lorden_delay_def} of the min-CuSum is achieved when the change occurs from the beginning, i.e.,   for every $i \in \mathcal{I}$  and $b>0$, 
 \begin{align} \label{CUSUM_worst case}
  \cJ_{i} \left[ \sigma(b) \right] = \Exp_i[\sigma(b)]
 \end{align} 
and  that~$\sigma (b_\alpha)$ minimizes~$\cJ_i$ in~$\cC (\alpha)$, for every~${i \in \mathcal{I}}$, to a first-order asymptotic approximation as~$\alpha \to 0$\@. In particular,  
 \begin{align} \label{CUSUM_AO}
  \cJ_{i} \left[ \sigma(b_\alpha) \right]  \sim  
  \frac{|\log \alpha|}{I_{i}}  \sim  \inf_{(T, D)  \in \mathcal{C}(\alpha)} \mathcal{J}_i[T]
 \end{align} 
  for every $i \in \mathcal{I}$  as $\alpha \to 0$\@.\\

\subsection{Min-CuSum for Sequential Change Diagnosis}
\label{subsec:intuition}

For any $b>0$, the stopping time $\sigma(b)$  is associated with a natural identification rule, which is   to select the post-change alternative  that corresponds to the largest CuSum  statistic at the time of stopping, i.e., 
\begin{equation}\label{CUSUM_identification_rule}
 \widehat{\sigma}(b)  \in    \argmax_{i \in \mathcal{I}} Y_i(\sigma(b)),
\end{equation}
solving ties, if any, in some arbitrary way. In other words, the post-change distribution  is declared to be the one that corresponds to the CuSum statistic that first crosses the threshold~$b$. Thus, for any $b>0$, $(\sigma(b), \widehat{\sigma}(b))$ is a sequential change diagnosis procedure, which we will still refer to as the min-CuSum, as its implementation is  the same as for the pure change-detection problem. The only (but main) difference is that in the context of change diagnosis the threshold~$b$ must be selected to satisfy a  false identification constraint  in addition to the false alarm constraint.\\

Before stating our theoretical results regarding the performance of min-CuSum as a change diagnosis procedure, it is useful to provide some intuition on the behavior of the CuSum statistics, $Y_i, i \in \cI$.



First of all, since each  $Y_i$ tends to return to 0 before the change, none of them will, in general, be large at the time of the change. Now, suppose that the true post-change density is $g = g_j$. Then, for every $n >\nu$  and $i \in \cI$ with $i \neq j$ we have 
\begin{align}
\begin{split}
\Exp_j[ \ell_j(n)] &= I_j   ,\\
 \Exp_j[ \ell_i(n)] &=   I_j -  I_{ji}. \\
 \end{split}
\end{align} 
This means that  the ``correct'' CuSum statistic,~$Y_j$,  grows  after the change essentially  like a random walk with positive drift, $I_j$. On the other hand, the post-change behavior of an ``incorrect'' CuSum statistic, $Y_i$,  depends on the sign of
$ I_j-  I_{ji}$, which reflects  whether the true post-change density, $g_j$,  is closer, in the sense of KL divergence, to the pre-change density,~$f$, or to the true post-change distribution,~$g_i$.  

\begin{itemize}
\item  If $ I_j < I_{ji}$, then  $g_j$ is closer to $f$ than $g_i$ and $Y_i$ behaves in the post-change regime as in the pre-change regime, i.e., it    tends to return to zero.




\item If $ I_j=  I_{ji}$, then $g_i$ and $f$ have the same ``distance'' from $g_j$ and  $Y_i$ becomes  a driftless random walk reflected at the origin. Therefore, while it still does not  systematically grow after the change,  it will  take much more time to return to 0 than before the change occurred.

\item
If $ I_j>  I_{ji}$, then $g_j$ is closer to $g_i$ than $f$ and $Y_i$ behaves similarly to $Y_j$, i.e., it grows essentially like a random walk with positive drift $I_j - I_{ji}$, which is notably less than the drift of $Y_j$. 
\end{itemize}

We continue by applying this discussion in the context of each of the scenarios presented in Section~\ref{subsec:multichannel} and Section~\ref{subsec:two-sided}, namely,  the single-fault and concurrent-fault multichannel problems, as well as the two-sided problem.



\begin{figure}[bt]
\begin{center}
\includegraphics[width=0.48\textwidth]{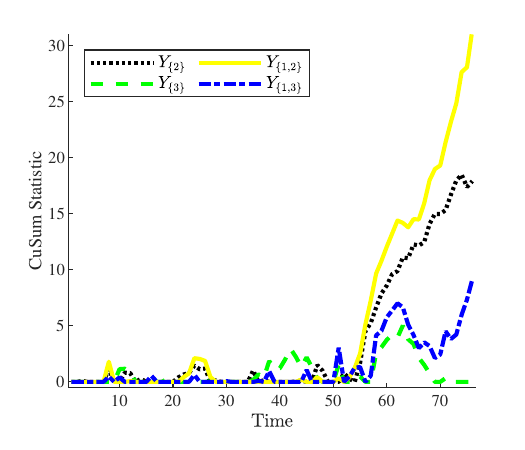}
\caption{Sample paths of some CuSum statistics in the concurrent-fault multichannel problem with $d=3$ channels and the change occurs in the first two of them, i.e.,  $g = g_{\{1,2\}}$, at time $\nu = 50$\@.}
\label{fig:intuition} 
\end{center}
\end{figure}

\begin{itemize}
\item We  consider, first, the \textit{single-fault} multichannel problem presented in Section~\ref{subsubsec:single_fault}, in which case $\cI=\{1, \ldots, d\}$, where~$d$ is the number of channels.  Then,  from~\eqref{multichannel KL} we can see that 
\begin{align} \label{ordering}
I_j &<  I_{ji} \quad \forall \; i, j \in \cI,\, i \neq j.
\end{align}
Thus, in this context, the values of all ``incorrect'' statistics will tend to be small  after the change. \\

\item We  consider, next, the \textit{concurrent-fault} multichannel problem presented in Section~\ref{subsubsec:multiple_fault}, with $d=3$ channels. Since this case is more complex, we assume for the sake of illustration that $p_i = p$ and~$q_i=q$ for $i = 1, \ldots, d, $ and that the pre- and post-change distributions are Gaussian with    standard deviation one and mean zero and one, respectively, i.e., \begin{equation}\label{eqn:Gaussian_distributions}
p  = \cN(0, 1) , \quad  q= \cN(1,1).
\end{equation}

We assume that the fault occurs in channels $1$ and $2$, i.e., $ g = g_{\{ 1, 2 \}}$. Then, we obtain from \eqref{eqn:KL_concurrent} and the symmetry of the Gaussian KL-number that  
\begin{align}
\begin{split}\label{eqn:relative_drifts}
I_{\{1,2\}} &< I_{\{1,2\}, \{3\}}  , \\
I_{\{1,2\}} &= I_{\{1,2\}, \{1,3\}}  ,\\
I_{\{1,2\}} &> I_{\{1,2\}, \{2\}}.\\
\end{split}
\end{align}
Thus,  we expect that, after the change,~$Y_{ \{3\}}$ will  tend to return  to   zero, whereas~$Y_{ \{1,3\} }$ will  not grow in a systematic way, but it will take much longer to return to zero\@. The statistic~$Y_{\{2\}}$  will grow roughly linearly but at a slower rate than the ``correct'' statistic, $Y_{\{1,2\}}$.  In Fig.~\ref{fig:intuition} we plot a path of these statistics when the change occurs at time~${\nu=50}$ with~${g =g_{\{1,2\}}}.$ \\

\item In the two-sided problem of Subsection~\ref{subsec:two-sided}, by~\eqref{eqn:two-sided_KL1}--\eqref{eqn:two-sided_KL2}, we have 
\begin{align}\label{eqn:KL_diff_two-sided}
\begin{split}
I_{ji} - I_j &= 
(\gamma_0 - \gamma_i)  \, \varphi'(\gamma_j) - ( \varphi(\gamma_0) - \varphi(\gamma_i)) \, \\
&= (\gamma_0 - \gamma_i) \left[  \varphi'(\gamma_j) - \frac{\varphi(\gamma_0) - \varphi(\gamma_i)}{\gamma_0 - \gamma_i} \right].
\end{split}
\end{align}
This implies that $I_j < I_{ji}$, which means that the ``incorrect" statistic tends to remain small after the change occurs.  To see that this inequality indeed holds, note that since $\varphi$ is  differentiable in the interior of $\Gamma$,  there exists some $\gamma_i^*$ between $\gamma_0$ and $ \gamma_i$ such that \begin{equation}\label{eqn:MVT}
 \varphi'(\gamma_i^*) = \frac{\varphi(\gamma_0) - \varphi(\gamma_i)}{\gamma_0 - \gamma_i},
\end{equation}
and, consequently, 
\begin{equation}
I_{ji} - I_j = (\gamma_0 - \gamma_i) \left[  \varphi'(\gamma_j) -  \varphi'(\gamma_i^*) \right].
\end{equation}
By the strict convexity of $ \varphi$ it follows that  if  $ \gamma_i < \gamma_0 < \gamma_j, $ then 
\begin{equation}
\varphi'(\gamma_j) -  \varphi'(\gamma_i^*)  > 0 \quad \& \quad  \gamma_0 - \gamma_i > 0  ,
\end{equation}
whereas if   $ \gamma_j < \gamma_0 < \gamma_i$, then 
\begin{equation}\label{eqn:two-sided-KL-final_eqn}
\varphi'(\gamma_j) -  \varphi'(\gamma_i^*) < 0 \quad \& \quad 
 \gamma_0 - \gamma_i < 0,
\end{equation}
which proves that, indeed, $I_j < I_{ji}$.
\end{itemize}





\section{Main Results}\label{sec:results}
In this section we establish our theoretical results.  We  first establish a non-asymptotic upper bound on the conditional probability  of mistakenly selecting $g_i$ when the true post-change density is $g_j$, given that there has not been a false alarm at the time of the change,
where $i, j \in \cI$ with $j \neq i$. This is established when   $g_j$  is closer---in the sense of Kullback-Leibler divergence---to the pre-change density $f$ than to $g_i$, i.e.,  ${I_j \leq I_{ji}}$. To obtain this result, we make the following second-moment assumption 
\begin{equation}\label{assum:KL_squared}
\int \left( \big( \log(g_i / f) \big)^+ \right)^2 g_i \, d\lambda < \infty,
\end{equation}
and we also make certain assumptions on the cumulant-generating function of~$\ell_i(n)$ under~$\Pro_{j}$,
 \begin{equation}
 \psi_{ij}(\theta) \equiv \log \left( \Exp_{j} \left[ \exp[\theta \, \ell_{i}(1)] \right] \right), \quad \theta \in \bR. \label{def:psi_ij}
 \end{equation}
 
The following theorem relates the relative distances between $g_i, g_j,$ and $f$ to the conditional probability of misidentication as a function of the threshold $b$.

 \begin{theorem}\label{thm:CUSUM_control}
Fix $i,  j \in \mathcal{I}, i \neq j$,  $b> 0$, and let~$\nu \geq 0$  such that \begin{equation}\label{condition on change points}
\Pro_{\infty}(Y_i(\nu) \geq x \; | \;  \sigma (b) > \nu) \leq  e^{-x} \quad \text{for all} \;  x \geq 0.
\end{equation} 

\begin{itemize}
\item[(i)] Suppose that~${I_j < I_{ji}}$ and that $\psi_{ij}$ has a positive root $r_{ij}$. 

\begin{itemize}
\item[(a)]
If~$r_{ij} > 1$, then 
\begin{align}
\label{eqn:exponential_bound1}
\begin{split}
  \Pro_{\nu, j}  ( \widehat{\sigma} (b)  = i \;   &| \; \sigma (b) > \nu )   \\
  &\leq  \frac{r_{ij}}{r_{ij} - 1}  \,  e^{-b}\,  (1+\phi(b)),
  \end{split}
\end{align}
 where $\phi(b)$ does not depend on $\nu$  and goes to $0$ as $b \to \infty$.

\item[(b)]If~$r_{ij} \in (0,1]$, then 
\begin{align} \label{eqn:exponential_bound2}
\begin{split}
  \Pro_{\nu, j} ( \widehat{\sigma} (b)  &= i \; | \; \sigma (b) > \nu ) \\
 & \leq ( r_{ij} + 1/I_j)  \, b\,  e^{-r_{ij}b} (1+\phi(b)),
  \end{split}
\end{align}
 where $\phi(b)$ is as above.\\
\end{itemize}
 
\item[(ii)]  If~${I_j = I_{ji}}$ and~$\psi_{ij}$ is finite on an open interval containing~$0$,  then 
  \begin{equation}\label{eqn:linear_bound} \Pro_{\nu, j}(\widehat{\sigma}(b) = i \; | \; \sigma(b) > \nu) \leq \widetilde{C}  \, b^{-1}\, \left(1 + \phi(b)\right) ,
  \end{equation}
   where~$\phi(b)$ is as in (i) and~$\widetilde{C}$ is a constant. 
\end{itemize}
\end{theorem}

\begin{IEEEproof}
Appendix~\ref{app:CUSUM_control}\@.
\end{IEEEproof}

To the best of our understanding, condition~\eqref{condition on change points} is assumed \text{implicitly} in the proof of the false identification control of the Vector CuSum~\cite[Theorem 2]{nikiforov2000simple}, as discussed  in~\cite{warner2022sequential}. 
This condition is  satisfied, trivially,  when  $\nu=0$\@.
   Indeed, it is well known that, for all $i \in \mathcal{I}$ and $ n \in \bN_0$, 
\begin{equation}\label{eqn:unconditional_bound}
\Pro_{\infty}(Y_i(n) \geq x) \leq e^{-x} \quad \text{ for all } x \geq 0.
\end{equation}
However, it is unclear to us  whether this is the case, in general, for other values of $\nu$.  We next show condition~\eqref{condition on change points} is indeed satisfied for every $\nu\geq 0$ in the multichannel setup in the case of a single fault. 


\subsection{The Single Fault Multichannel Setup}\label{sec:multichannel}

We next apply Theorem \ref{thm:CUSUM_control} to the case of the   \textit{single-fault} multichannel problem introduced in Section~\ref{subsubsec:single_fault}.
In this setup, we  show that the worst-case---with respect to the change-point---conditional probability of misidentification of the min-CuSum decays exponentially fast in the threshold $b$ and,  as a result,  the min-CuSum  achieves  \eqref{infimum}, simultaneously for every~${i \in \mathcal{I}}$, 
as $\alpha$ goes to 0 sufficiently faster than $\beta$.


\begin{corollary}\label{coro:multichannel_control}
Consider the  single-fault multichannel problem of Section~\ref{subsubsec:single_fault} and suppose that, for every $i \in \cI$, 
\begin{align*}
  \int_{\bS_i} \left( \log^2 (q_i/ p_i) \right)  q_i \, d \lambda_i< \infty. 
\end{align*}

\begin{enumerate}
\item[(i)] For any~$b > 0$,
\begin{align*}
\max_{j \in \mathcal{I}} \sup_{\nu \in \bN_0}  \Pro_{\nu,j} ( \widehat{\sigma}(b) \neq j\,  | \, \sigma(b) > \nu) 
\leq \\
C \, b \, e^{-b} \, (1 + \phi(b)),
\end{align*}
where $\phi(b)$ is a function of $b$ that does not depend on $\nu$ and  goes to $0$ as ${b \to \infty}$, and \begin{equation}
C \equiv ({|\mathcal{I}|}-1) \left( 1 + \max_{i \in \mathcal{I}}(1/I_i) \right).
\end{equation}

\item[(ii)] Suppose also that~$\alpha$ and~$\beta$ satisfy 
\begin{equation}\label{CUSUM beta condition}
C  \, b_{\alpha} \, e^{-b_{\alpha}} \left( 1 + \phi
(b_{\alpha}) \right) \leq \beta,
\end{equation}
where~$C$ and~$\phi$ are defined as above, and~$b_{\alpha}$ is defined as in~\eqref{eqn:b_alpha}.  Then, \begin{equation}
\left( \sigma(b_{\alpha}), \widehat{\sigma}(b_{\alpha}) \right) \in \cC(\alpha, \beta) \quad \forall \, \alpha \in (0,1)
\end{equation}
and, for every $i \in \mathcal{I}$, \begin{equation}
\cJ_i[\sigma(b_{\alpha})] \sim \inf_{(T,D) \in \cC(\alpha, \beta)} \cJ_i[T] \quad \text{ as } \alpha \to 0.
\end{equation}

\end{enumerate}

\end{corollary}

\begin{IEEEproof}
(i) By Boole's inequality, for any $j \in \cI$ we have 
\begin{equation*}
\Pro_{\nu, j} \left( \widehat{\sigma}(b) \neq j \, | \, \sigma(b) > \nu \right) \leq \sum_{i \neq j} \Pro_{\nu, j} \left( \widehat{\sigma}(b) = i \, | \, \sigma(b) > \nu \right).
\end{equation*}
Therefore,  it suffices to show that, for every $j \in \cI$, 
\eqref{eqn:exponential_bound2} holds with $r_{ij}=1$ for every   $ \nu \in \bN_0$, $i \in \mathcal{I}$, $ b > 0$. In what follows we fix $\nu, i,j, b$.

 We   first  show that condition~\eqref{condition on change points} holds for every $\nu\geq 0$. Indeed, for all   $x \geq 0$ we have 
\begin{align*}
& \Pro_{\infty}(Y_i(\nu) \geq x  \, | \, \sigma (b) > \nu) \\
 &= \Pro_{\infty}(Y_i(\nu) \geq x \, |\,  \sigma_j(b) > \nu \;  \text{for every} \; j \in \mathcal{I}) \\
 &= \Pro_{\infty}(Y_i(\nu) \geq x \; | \; \sigma_i (b) > \nu)  \\
 &\leq \Pro_{\infty}(Y_i(\nu) \geq x).
 \end{align*}
The  first equality follows by the definition of $\sigma(b)$ and the second  by the independence of the channels, whereas the inequality follows by \cite[Theorem 1]{pollak1986convergence}, since  $Y_i$
is a non-negative stochastically monotone Markov process (Lemma~\ref{lemma:Multichannel_assumption} in Appendix~\ref{app:lemmas}). 

We have already seen that in this multichannel setup we have $I_{j}<I_{ji}$ for every $i, j \in \mathcal{I}$ with $i \neq j$ (recall \eqref{ordering}). Therefore, in view of Theorem~\ref{thm:CUSUM_control}(i), it remains to show that each $\psi_{ij}$ has root equal to~$1$.  Indeed, for each $i \in \mathcal{I}$, the  likelihood ratio $\Lambda_i(n)$ defined in~\eqref{LLR}  takes the following form:
\begin{align*}
\Lambda_i(n) &= \frac{q_i(X_{i,n})}{p_i(X_{i,n})}, \quad  \; n \in \bN, 
\end{align*}
and, consequently, 
\begin{align*}
\exp [\psi_{ij}(1)] &= \Exp_j[\exp [\ell_i(1)]] \\
&=  
\Exp_j \left[  \frac{q_i(X_{i,n})}{p_i(X_{i,n})} \right] =   \int (q_i/ p_i) \, p_i \, d \lambda_i= 1.
\end{align*}
 
(ii)     For any $\alpha, \beta \in (0,1)$, we have 
$ \cC(\alpha, \beta) \subseteq \cC(\alpha)$,
 therefore
\begin{equation} 
 \inf_{(T, D)  \in \cC(\alpha, \beta, \bN_0)}   \cJ_i[T] \geq 
  \inf_{(T,D)  \in \cC(\alpha)}   \cJ_i[T].
 \end{equation} 
In view of~\eqref{CUSUM_AO}, it suffices to show that  ${\sigma(b _\alpha) \in \cC(\alpha, \beta)}$ when~$\beta$ satisfies~\eqref{CUSUM beta condition}, which follows by~Corollary~\ref{coro:multichannel_control}(i)\@.
\end{IEEEproof}



\subsection{The Two-Sided Setup}
We next apply Theorem  \ref{thm:CUSUM_control} to the  two-sided setup of Subsection~\ref{subsec:two-sided}. In this case, we have already seen that condition~\eqref{ordering} holds for every $i, j \in \mathcal{I}$  with  $i \neq j$. Moreover, by the definition of $\psi_{ij}$ in \eqref{def:psi_ij} we have 

\begin{align} \label{eqn:exp_fam_cumulant}
\begin{split}
\psi_{ij}(\theta) &\equiv \log \int \exp \left\{ \theta  \log \frac{h_{\gamma_j}}{f}  \right\} \; h_{\gamma_i} \; \nu(dx) \\
&= \varphi(\theta(\gamma_j - \gamma_0) + \gamma_i)  -\theta(\varphi(\gamma_j) - \varphi(\gamma_0)) - \varphi(\gamma_i) ,
\end{split}
\end{align}

Since $\psi_{ij}$ is finite,  strictly convex, and differentiable on $\Gamma$ with 
$${\psi_{ij}'(0) = \Exp_j[\ell_i(1)]} =  I_1 - I_{12} <0  ,
$$
 $\psi_{ij}$ has a positive root when  $\Gamma=\bR$. 
  Thus, the probability of misidentification will decay exponentially in the threshold~$b$, depending on the roots of~$\psi_{ij}.$ In the following corollary, we specialize this result to the mean-shift of a Gaussian sequence,  where we explicitly compute the misidentification rate from Theorem~\ref{thm:CUSUM_control}\@. 

\begin{corollary} Consider the two-sided problem with $h_{\gamma} = N(\gamma, 1)$ and $ \gamma_0 = 0, $ and  $\gamma_1 < 0 < \gamma_2$.  For any $ b > 0$  and $\nu \in \bN_0$ for which  \eqref{condition on change points} holds  we have 
\begin{equation}
\max_{j = 1, 2} \sup_{\nu \in \mathbb{N}_0} \Pro_{\nu, j}\left( \widehat{\sigma}(b) \neq j \, | \, \sigma(b) > \nu \right) \leq C^*  e^{-b} ( 1 + \phi(b) )  ,
\end{equation}
where $$ C^* = 1 + \frac{1}{2}\max \left\lbrace \frac{|\gamma_1|}{|\gamma_2|}, \frac{|\gamma_2|}{|\gamma_1|}\right\rbrace $$
and~$\phi(b)$ is as in Theorem~\ref{thm:CUSUM_control}\@.  
\end{corollary}
\begin{IEEEproof}
It follows from \eqref{eqn:KL_diff_two-sided}--\eqref{eqn:two-sided-KL-final_eqn}  that~\eqref{ordering} holds, so it remains to compute the roots of~$\psi_{ij}$ using~\eqref{eqn:exp_fam_cumulant}. Since~${\varphi(\gamma) = \gamma^2/2},$ solving $\psi_{ij}(\theta) = 0$ for~$\theta$ yields~${\theta = 1 + 2 |\gamma_j / \gamma_i|}$. Finally, using Theorem~\ref{thm:CUSUM_control}(i)(a), the proof is complete. 
\end{IEEEproof}

\section{Simulation Studies}\label{sec:simulation}

In this section, we simulate two variations of the  multichannel problem introduced in Section~\ref{subsec:multichannel}, with three data channels, i.e.,~$d=3$\@. In the first study we consider the single-fault multichannel problem of Section~\ref{subsubsec:single_fault}, i.e., $\mathcal{I} = \{1, \ldots, d \}$, and in the second study  the concurrent-fault multichannel problem of Section~\ref{subsubsec:multiple_fault}, i.e.,~$\mathcal{I}$ is the set of all non-empty subsets of~$ \{ 1, \ldots, d \}$. In both studies,  the pre- and post-change distributions of the channels are assumed to be Gaussian, as in~\eqref{eqn:Gaussian_distributions}.

\subsection{Single-Fault Multichannel Problem}
\begin{figure}[bt]
\begin{center}
\includegraphics[width=0.48\textwidth]{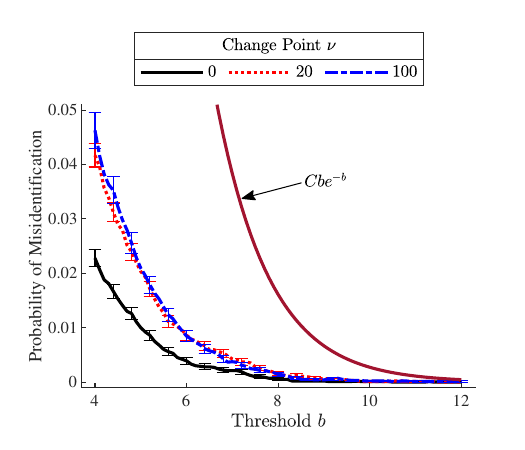}
\caption{For the single-fault multichannel problem, the conditional probability of misidentification, given that there was not a false alarm, is plotted as a function of the threshold~$b$, for several values of the change point, ${\nu = 0, 20, 100}$\@. Standard errors are represented with error bars, and only every fourth error bar is included for the sake of visual clarity. The first-order term of the theoretical upper bound for the worst-case conditional probability of misidentification is also plotted.}
\label{fig:study2} 
\end{center}
\end{figure}

In this study, we consider the case where the fault occurs in channel 1 only i.e.,~
$ g= g_1$.  Due to symmetry, there is no reason to consider the case where  $g=g_2.$  For~$\nu = 0, 20$, and~$100$, we simulate~$10^4$ paths of the CuSum statistics~$Y_i, i \in \cI$, where $ g= g_1,$ in order to estimate
\begin{equation}
\Pro_{\nu,   1 } \left( \widehat{\sigma}(b) \neq  1  \, | \, \sigma(b) > \nu \right)
\end{equation}
as a function of~$b$. We report the results in Fig.~\ref{fig:study2}\@. For each change point, the conditional probability of misidentification goes to zero as $b$ increases, and the decay appears exponential in~$b$. The conditional probability of misidentification increases as $ \nu$ increases, with a noticeable difference between the cases of ${\nu = 0}$ and ${\nu = 20}$\@. However, the estimated probabilities when ${\nu = 20}$ and ${\nu = 100}$ are nearly indistinguishable. We also include a curve obtained from the first-order term~${C} be^{-b}$ in the upper bound from Corollary~\ref{coro:multichannel_control}, where in particular straight-forward calculation shows that ${C = 6}$\@. 

\subsection{Concurrent-Fault Multichannel Problem}
In this study, we analyze the case where the fault occurs in channels 1 and 2 simultaneously, i.e.,~$g = g_{\{ 1,2\} }$, in the spirit of Section~\ref{subsec:intuition}\@. For~${\nu = 0, 20,} $ and~$100$, we simulate~$10^4$ paths of the CuSum statistics $Y_i, i \in \mathcal{I},$ in order to estimate the overall probability of misidentification, 
\begin{equation}\label{total_probs}
    \Pro_{\nu, \{ 1,2 \}} \left( \widehat{\sigma}(b) \neq \{1, 2\} \, | \, \sigma(b) > \nu \right)  ,
\end{equation}
as well as the probability of incorrectly identifying each of~$g_{\{2\}}, g_{\{1,3\}},$ and $g_{\{3\}}$ as the post-change distribution, i.e., 
 \begin{equation}\label{marginal_probs}
\Pro_{\nu, \{ 1,2 \}} \left( \widehat{\sigma}(b) = k \, | \, \sigma(b) > \nu \right) 
\end{equation}
for~${k = \{ 2 \}, \{ 1,3 \} ,  \{ 3 \}}$.

\begin{figure}[bt]
\begin{center}
\includegraphics[width=0.48\textwidth]{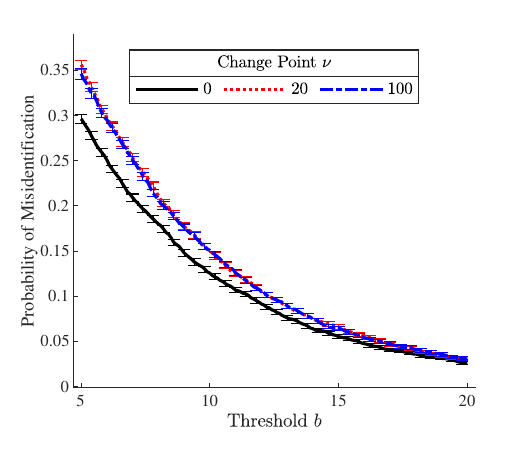}
\caption{For the concurrent-fault multichannel problem, the overall conditional probability of misidentification is plotted, when the true post-change distribution is $g_{\{1,2 \}}$. The curves are plotted as functions of the threshold~$b$, and the three curves correspond to the cases of $\nu = 0, 20,$ and $100$\@. Standard errors are represented with error bars, and only every fourth error bar is included for the sake of visual clarity.}
\label{fig:study1A} 
\end{center}
\end{figure}

In Fig.~\ref{fig:study1A}, we report the overall probability of misidentification \eqref{total_probs}
 for $ \nu = 0, 20,$ and~$100$\@. As in the single-fault study, this is smaller for $\nu = 0$ than for $\nu = 20$, whereas it is essentially the same for $ \nu =20$ and $ \nu = 100$. 

\begin{figure}[bt]
\begin{center}
\includegraphics[width=0.48\textwidth]{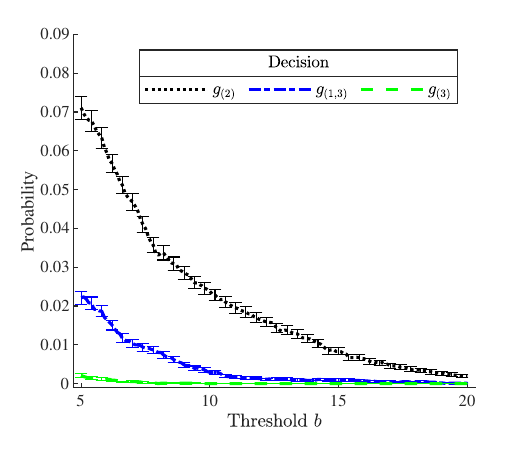}
\caption{The three curves correspond to the conditional probabilities of identifying $g_{\{2\}}, g_{\{1,3\}}$, and $g_{\{3\}}$ as the post-change distribution, given that there was not a false alarm, when the true post-change distribution is $g_{\{1,2 \}}.$ These curves are plotted as functions of the threshold~$b$, and the change time is~$\nu =100$\@. 
Standard errors are represented with error bars, and only every fourth error bar is included for the sake of visual clarity.}
\label{fig:study1B} 
\end{center}
\end{figure}

In Fig.~\ref{fig:study1B}, we report the ``partial'' probabilities of misidentification in ~\eqref{marginal_probs} for ~${\nu =100}$\@. We see that  these are ordered in the way suggested by~\eqref{eqn:relative_drifts},
with~$g_{ \{ 2 \} }$ being more frequently incorrectly identified as the post-change distribution than~$g_{ \{ 1,3 \} }$, and the latter more frequently identified as the post-change distribution than~$g_{ \{ 3 \} }$.




\section{Conclusion}\label{sec:conclusion}
In this work, we  study  the min-CuSum algorithm in the context of the sequential change diagnosis, where it is desirable to not only detect the change but also to isolate the post-change regime among a finite set of alternatives. In the special case that the data are collected over independent channels and the change occurs in one of the channels, we show that the conditional probability of misidentification decays exponentially in the threshold, \textit{uniformly for all possible change points}.  Our results imply a similar result for the two-sided problem in an exponential family, under a certain technical condition whose verification remains an open problem. 



The sequential change diagnosis algorithm we consider in this work can be naturally extended  to more general, pre-change and post-change models that  allow for composite hypotheses and/or non-iid data.  A theoretical analysis in such cases 
is an interesting subject for future work. Moreover, another interesting direction for further research is  the study of  \textit{decentralized} sequential change diagnosis algorithms in the context of the multichannel problem. For example,
suppose it is known a priori that the change will occur in  exactly $M$ channels, and  consider the algorithm that  stops as soon as $M$ CuSum statistics have crossed the threshold, as in~\cite{banerjee2016decentralized}. Our results in the present paper cover the special case where $M=1$. An analysis in the case where the change can occur in  more than one channels simultaneously, i.e., $M>1$, is an open problem.


 


\appendices

\section{Properties of CuSum stopping times}
\label{app:cusums}

 Throughout this appendix we fix~$b>0$ and $i, j \in \mathcal{I}$ such that $i \neq j$. We recall from Section~\ref{sec:CUSUM} that $\sigma_{i}(b)$ is the CuSum stopping time,  with threshold $b$,   for detecting a change from $f$ to $g_i$, i.e.,  
  \begin{align*}
 \sigma_{i}(b) &\equiv  \inf \lbrace n \in \bN : Y_{i}(n) \geq b  \rbrace .
\end{align*} We use the following notation for the  log-likelihood ratios for $n \in \mathbb{N}$: 
 \begin{equation}
 Z_i(n) \equiv \sum_{k=1}^n \ell_i(k).
 \end{equation}
It is known (e.g.,~\cite[Section 8.2.6.5]{tartakovsky2014sequential}) that  the expected value of $\sigma_i$  admits the following lower bounds when the change occurs at at $\nu=0$ and the true post-change hypothesis is $g_j$:
\begin{itemize}
\item If~$I_j < I_{ji}$ and $\psi_{ij}$ has a positive root $r_{ij}$, then \begin{equation}\label{eqn:exponential_lower_bound}
\Exp_j[\sigma_{i}(b)] \geq e^{r_{ij}b}.
\end{equation}

\item If $I_j = I_{ji}$ and $\psi_{ij}$ is finite on an interval containing zero, then 
\begin{equation}\label{eqn:quadratic_lower_bound}
 \Exp_j[\sigma_i(b)]  \geq \frac{b^2}{\Exp_j[\ell_i^2(1)]}. 
\end{equation}
\end{itemize}


We next develop analogous lower bounds in the case that  $Y_i$  is initialized from some value~$x$ in $[0,b]$. Specifically, we denote by $L_{ij}(x;b)$ the expectation under $\Pro_j$ of  $\sigma_i(b)$, when  $ Y_{i}$ is initialized from  some $x \in [0, b]$, i.e., 
\begin{equation} \label{L}
L_{ij}(x;b) \equiv \Exp_j[\sigma_{i}(b) \, |\,  Y_{i}(0) = x], \quad x \in [0, b].  
\end{equation}
Clearly,   $L( \, \cdot \, ;b)$ is a non-increasing function (see also Lemma~\ref{lemma:Multichannel_assumption}), so that 
\begin{equation}
L_{ij}(x;b) \leq L_{ij}(0;b) \equiv \Exp_j[\sigma_{i}(b)], \quad x \in [0,b].\\ 
 \end{equation}
To  provide explicit lower bounds on the function~$L_{ij}(x;b)$, as well as to bound $\Exp_j[\sigma_j(b)]$, we introduce the following quantities:
\begin{align}\label{def: omegas}
\begin{split}
B_{j} &\equiv \Exp_j\Big[\left( (\ell_j)^+ \right)^2\Big]/I_j^2  , \\
\omega_{ij} &\equiv \sup_{t \geq 0} \Exp_j \left[ \ell_{i}(1) - t \, | \, \ell_{i}(1) \geq t \right] , \\
\tilde{\omega}_{ij} &\equiv -  \inf_{t \leq 0} \Exp_j \left[ \ell_{i}(1) - t \, | \, \ell_i(1) - t \leq 0 \right]. 
\end{split}
\end{align}
The quantity $B_{j} $ is Lorden's upper bound on the excess over a positive threshold  of the random walk $Z_{j}$ under~$\Pro_j$ \cite{lorden1970excess}. In particular, we have 
\begin{equation}\label{eqn:cusum_delay_upper_bound}
\Exp_j[\sigma_j(b)] \leq \frac{b}{I_{j}} + B_{j}.
\end{equation}

The quantity $\omega_{ij}$  is a bound for the expected overshoot of the random walk $Z_{i}$ under~$\Pro_j$ above a non-negative threshold, and it is finite when $ \psi_{ij}$ has a positive root (see~\cite{bryson1974heavy}).  Similarly, $\tilde{\omega}_{ij}$ is a bound on the expected undershoot  $ \omega_{ij}$ of the random walk $Z_{i}$ above a non-positive threshold under~$\Pro_j$, and it is finite when  $\psi_{ij}$ is finite on an open interval containing zero.



\begin{lemma}\label{L_Lemma}
\mbox{}  
\begin{itemize}
\item[(i)] If~$I_j < I_{ji}$ and $\psi_{ij}$ has a positive root $r_{ij}$, then 
$$ L_{ij}(x ; b) \geq \l_{ij}(x; b), \quad \text{ for all } x \in [0, b] , $$ where
\begin{align} \label{def: l}
\begin{split}
l_{ij}(x; b)  &\equiv  \frac{x - e^{-r_{ij}(b- x)}(b+ \omega_{ij})}{I_{ji} - I_j}  \\
& + (1 - e^{-r_{ij}(b- x)}) \; L_{ij}(0; b).
\end{split}
\end{align}

\item[(ii)] If~${I_{ji} = I_j}$ and~$\psi_{ij}$ is finite on an interval containing zero, then
\begin{equation}
 L_{ij}(x \, ; \, b) \geq u_{ij}(x \, ; \, b) ,
 \end{equation}
 where
\begin{align}\label{def: u}
 u_{ij}(x ;  b)  \equiv \frac{b - x}{b + \omega_{ij} + \tilde{\omega}_{ij}} \times \left( \frac{xb}{\Exp_j[\ell_i^2(1)]} + L_{ij}(0;b) \right).
 \end{align}
\end{itemize}
\end{lemma}

\begin{IEEEproof}
We denote by~${(T_i(x), \delta_{i}(x))}$  the stopping time and decision of the sequential probability ratio test~(SPRT) of~$g_i$ vs~$f$ with boundaries zero and $b$ when the corresponding statistic is initialized from~${x \in [0,b]}$, i.e.,  
\begin{align*}
T_{i}(x) &= \inf \lbrace n \in \bN : Z_{i}(n) \notin (-x,  b- x)   \rbrace  \\
\delta_{i}(x) &= \begin{cases} 
0, \quad \text{if} \quad  Z_{i}(T_{i}(x)) \leq -x  \\
1,  \quad \text{if} \quad Z_{i}(T_{i}(x)) \geq b - x, \\
\end{cases}
\end{align*}
and we set 
$$ p(x) \equiv \Pro_j(\delta_{i}(x) = 1).$$


After initializing the CuSum from $x $, the CuSum statistic is equivalent to an SPRT with stopping boundaries being the positive threshold and zero. If the CuSum statistic crosses the positive threshold, the process terminates, whereas if it reaches zero, the process ``restarts", and as a result \cite[Section 5.2.2.1]{basseville1993detection} we have:
\begin{equation}\label{eqn:ARL_decomposition}
 L_{ij}(x ; b) = \Exp_j \big[ T_{i}(x) \big] + (1- p(x)) L_{ij}(0 ; b), \quad x \in [0, b]. 
\end{equation}

We now prove parts (i) and (ii) separately, in each case providing first an upper bound on $p(x)$ and then a lower bound on $\Exp_j[T_i(x)].$\\

\textbf{Proof of (i)} This proof is similar to parts of the proof of~\cite[Theorem 2]{nikiforov2000simple}. Since 
$$ \left\lbrace \exp \left[ r_{ij} \, Z_{i}(n)\right], n \in \mathbb{N} \right\rbrace $$ is a non-negative $\Pro_j$-martingale with expectation $1$,  by the definition of $p(x)$ and  Ville's supermartingale inequality we have 
\begin{equation}\label{process_bound}
p(x) \leq  \Pro_j \left( \sup_{n \geq 1} Z_{i}(n) > b - x \right) \leq e^{-r_{ij}(b-x)}.
\end{equation} Comparing with the desired inequality, it is enough to show that \begin{equation}
\Exp_j[T_{i}(x)] \geq \frac{x -  p(x)(b + \omega_{ij})}{I_{ji} - I_j}, \quad \forall x \in [0,b]. 
\end{equation} This can be shown by modifying the proof of the upper bound for~$\Exp_j \big[ Z_i( T_{i}(x))  \big]$ in \cite[Section 3.1.2]{tartakovsky2014sequential} in view of the fact that  $\Exp_j[\ell_i(1)] = I_{j} - I_{ji} < 0.$ \\



\textbf{Proof of (ii)}
Since $ \Exp_j[\ell_i(1)] = 0$ and $\Exp_j[T_i(x)]<\infty$, 
by Wald's identity and the law of total probability we have   \begin{align}\label{eqn:conditional}
\begin{split}
0&=\Exp_j[Z_i(T_i(x))]  \\
&=\Exp_j[Z_i(T_i(x)) \, | \, Z_i(T_i(x)) \geq b-x] \; p(x) \\   &\quad +\Exp_j[Z_i(T_i(x)) \, | \, Z_i(T_i(x)) \leq -x] \;  (1 - p(x)).
\end{split}
\end{align}  
By the definition of conditional expectation, 
\begin{align}\label{eqn:conditional_bound}
\begin{split}
\Exp_j[Z_i(T_i(x)) \,|\,Z_i(T_i(x)) \geq b - x] &\geq b - x \\
\Exp_j[Z_i(T_i(x)) \,|\,Z_i(T_i(x)) \leq - x] &\leq  -x .
\end{split}
\end{align}
By the definitions of $ \omega_{ij}$ and $ \tilde{\omega}_{ij}$ given in~\eqref{def: omegas} we also have 
\begin{align}\label{eqn:conditional_bound2}
\begin{split}
\Exp_j[Z_i(T_i(x)) \,|\,Z_i(T_i(x)) \geq b - x] &\leq b - x + \omega_{ij}\\
\Exp_j[Z_i(T_i(x)) \,|\,Z_i(T_i(x)) \leq - x] &\geq - (x + \tilde{\omega}_{ij}).
\end{split}
\end{align}
Solving~\eqref{eqn:conditional} for $p(x)$ and using~\eqref{eqn:conditional_bound} and~\eqref{eqn:conditional_bound2}, we conclude that 
\begin{equation}\label{eqn:prob2_lowerbound}
p(x) \geq \frac{x}{b + \omega_{ij} + \tilde{\omega}}_{ij}. 
\end{equation}
By solving~\eqref{eqn:conditional} for $ 1 - p(x)$, we can use~\eqref{eqn:conditional_bound} and~\eqref{eqn:conditional_bound2} to show that
\begin{equation}\label{eqn:prob1_lowerbound}
1 - p(x) \geq \frac{b-x}{b + \omega_{ij} + \tilde{\omega}_{ij}} .
\end{equation} 
Combining~\eqref{eqn:prob1_lowerbound} with~\eqref{eqn:ARL_decomposition} we have 
\begin{equation}\label{eqn: halfway_to_lowerbound}
 L_{ij}(x;b) \geq \Exp_j[T_i(x)] +  \frac{b-x}{b + \omega_{ij} + \tilde{\omega}_{ij}}  \, L_{ij}(0;b).
 \end{equation} It remains to lower bound $ \Exp_j[T_i(x)]$.
Since, by assumption,  $\Exp_j[\ell_i^2(1)]$ is finite,  by Wald's identity for the second moment
  we obtain
\begin{align}\label{zero_lb}
\begin{split}
& \Exp_j[T_i(x)] = \frac{\Exp_j[Z_i^2(T_i(x))]}{\Exp_j [\ell_i^2(1)]} \\
& = \frac{\Exp_j[Z_i^2(T_i(x)) \,|\, Z_i(T_i(x)) \geq b -x] \, p(x)}{\Exp_j [\ell_i^2(1)]} \\
&  \quad + \frac{\Exp_j[Z_i^2(T_i(x)) \,|\, Z_i(T_i(x)) \leq -x] \, (1 - p(x))}{\Exp_j [\ell_i^2(1)]} \\
& \geq \frac{(b - x)^2\,  p(x) + x^2\,  (1 - p(x))}{\Exp_j [\ell_i^2(1)]},
\end{split}
\end{align}
where for the inequality we use~\eqref{eqn:conditional_bound}. 
Combining~\eqref{zero_lb} with 
\eqref{eqn:prob2_lowerbound}--\eqref{eqn:prob1_lowerbound}, we obtain the desired  inequality.
\end{IEEEproof}

By direct calculation, we can obtain the following bounds on the derivatives of $u_{ij}$ and $l_{ij}$, will will  be useful later. 

\begin{lemma}  \label{lemma:derivatives}
\begin{itemize}
\item If $I_{j} < I_{ji}$ and $ \psi_{ij}(\theta)$ has positive root $r_{ij}$, then 
  \begin{align*}
0 <   -l'_{ij}(x;b) &\leq r_{ij} e^{-r_{ij}(b - x)} \left( \frac{b + \omega_{ij}}{I_{ji} - I_j} + L_{ij}(0;b) \right).
\end{align*}  

   \item If $I_j = I_{ij}$ and $\psi_{ij}$ is finite on an open interval containing zero, then 
  \begin{align*}
- u_{ij}'(x\,;\,b) 
&= \frac{L_{ij}(0;b)}{ b + \omega_{ij} + \tilde{\omega}_{ij}} \,  + \,  \frac{xb - (b^2 - xb)}{\Exp_j[\ell_i(1)^2] \, ( b + \omega_{ij} + \tilde{\omega}_{ij})} \\
& \leq \frac{L_{ij}(0;b)}{b + \omega_{ij} + \tilde{\omega}_{ij}}  +\frac{xb}{\Exp_j[\ell_i(1)^2] \, (b + \omega_{ij} + \tilde{\omega}_{ij})}.
\end{align*}
\end{itemize}
\end{lemma}


\section{Proof of Theorem~\ref{thm:CUSUM_control}}\label{app:CUSUM_control}

 We start with a lemma that shows that  when $b$ is large enough, a relevant conditional probability is non-zero. 

\begin{lemma}\label{lem:nonzero_prob}
Fix $i, j\in \mathcal{I}$  with $i \neq j.$
Assume either that~$I_{ji} > I_j$ and $\psi_{ij}$ has a positive root, or that~${I_{ji} = I_j}$ and~$\psi_{ij}$ is finite around zero. Let~$\nu$ be such that~\eqref{condition on change points} holds, at least for~$b$ large enough. 
If~$b$ is large enough,  
then 
$$ \Pro_{\nu, j}(\sigma_i(b) > \sigma_j(b) \, | \, \sigma(b) > \nu) > 0. $$ 
\end{lemma}

\begin{IEEEproof}
To make the notation lighter, we write $\sigma$ and $\sigma_i$  in place of~$\sigma(b)$ and $\sigma_i(b)$.
We first observe that, by the law of total probability, 
\begin{align*}
\Pro_{\nu, j}(\sigma_i > \sigma_j \, | \, \sigma > \nu)  & \geq \Pro_{\nu, j}( Y_i(\nu) \leq b/2 \, | \,\sigma > \nu) \,\\
& \quad \times 
\Pro_{\nu, j}(\sigma_i > \sigma_j \, | \, Y_i(\nu) \leq b/2, \sigma > \nu)  , 
\end{align*}
so it is enough to show that each of the two terms in the lower bound is non-zero. Since $Y_i(\nu)$ has the same distribution under $ \Pro_{\infty}$ as under $ \Pro_{\nu, j}$, condition \eqref{condition on change points}  implies that  \begin{equation}
\Pro_{\nu, j}( Y_i(\nu) \leq b/2 \, | \,\sigma > \nu) \geq 1 - e^{-b/2}>0.
\end{equation}
Thus, it remains to show that the second term in the upper bound is also positive.  For this, 
we  apply the law of total expectation and we write 
\begin{align}\label{0}
\begin{split}
\Pro_{\nu, j}(\sigma_i > \sigma_j \,  | \, Y_i(\nu) \leq b/2, \sigma > \nu) = \\
\Exp_{\nu, j}\left[ \Pro_{\nu, j}(\sigma_i > \sigma_j \, | \, \cF_{\nu} ,\sigma > \nu ) \, \big| \, Y_i(\nu) \leq b/2, \, \sigma > \nu \right].
\end{split}
\end{align}
When $\sigma > \nu$, the stopping times~$\sigma_i, \sigma_j$ depend on~$\cF_{\nu}$ only through~$Y_i(\nu)$ and~$ Y_j(\nu)$, thus we can write 
\begin{align}\label{1}
\begin{split}
\Exp_{\nu, j}\left[ \Pro_{\nu, j}(\sigma_i > \sigma_j \, | \, \cF_{\nu}, \sigma > \nu ) \, \big| \, Y_i(\nu) \leq b/2, \, \sigma > \nu \right]  \\
=\Exp_{\nu, j}\left[ G(Y_i(\nu), Y_j(\nu) ) \, \big| \, Y_i(\nu) \leq b/2, \, \sigma > \nu \right],
\end{split}
\end{align}
where $G$ is defined as the function 
\begin{align*}
G(y_i, y_j) &\equiv \Pro_{\nu, j}(\sigma_i > \sigma_j \, | \, Y_i(\nu) = y_i, \, Y_j(\nu) = y_j ) \\
&= \Pro_j(\sigma_i > \sigma_j \, | \, Y_i(0) = y_i, \, Y_j(0) = y_j ) ,
\end{align*}
where the equality is due to the data being iid after the change. 
This function,~$G$, is decreasing in~$y_i$ and increasing in~$y_j$. Thus, for all values of $Y_i(\nu), Y_j(\nu)$ in ${\lbrace Y_i(\nu) \leq b/2, \, \sigma > \nu \rbrace}$, 
\begin{align}\label{2}
 G(Y_i(\nu), Y_j(\nu)) &\geq G(b/2, 0).
 \end{align}
In particular, \begin{equation}\label{3}
    G(b/2, 0) = \Pro_j(\sigma_i > \sigma_j \, | \, Y_i(0) = b/2, \, Y_j(0) = 0).
\end{equation}
Combining~\eqref{0},~\eqref{1},~\eqref{2}, and~\eqref{3}, and the monotonicity of expectation, 
it follows that 
\begin{align*}
&\Pro_{\nu, j}(\sigma_i > \sigma_j | \, Y_i(\nu) \leq b/2, \sigma > \nu) \\
&\geq \Pro_j(\sigma_i > \sigma_j \, | \, Y_i(0) = b/2) ,
\end{align*}
where it is understood that $Y_j(0) = 0$\@.
Now, we claim that \begin{equation*}
 \Pro_j(\sigma_i(b) > \sigma_j(b) \, | \, Y_i(0) = b/2) \geq \Pro_j(\sigma_i(b/2) > \sigma_j(b)).
\end{equation*}
Indeed, $\sigma_i(b)$ with $Y_i(0) = b/2 $ is stochastically larger than~$ \sigma_i(b/2)$ with initialization at zero since, in both cases, the statistic $Y_i$ must increase by~$b/2 $ to reach the threshold. However, in the first case, the statistic~$Y_i$ can decrease below its initialization at~$b/2$, whereas~$Y_i$ is reflected at its initialization in the second case, and is therefore pathwise larger. 

Thus, it remains to show that  \begin{equation*}
\Pro_j(\sigma_i(b/2) > \sigma_j(b)) > 0. 
\end{equation*}
In order to do so,  we argue by contradiction and  suppose that  $\sigma_i(b/2) \leq  \sigma_j(b)$ with probability one under~$\Pro_j$. It follows from this assumption and~\eqref{eqn:cusum_delay_upper_bound} that 
\begin{equation} \label{contr}
 \Exp_j[\sigma_i(b/2)] \leq 
\Exp_j[\sigma_j(b)] \leq \frac{b}{I_j} + B_j.
\end{equation}
However, by either~\eqref{eqn:exponential_lower_bound} or~\eqref{eqn:quadratic_lower_bound}, depending on whether~$I_{ij}$ is greater than or equal to~$I_j$, we have that 
\begin{align*}
 \text{ either} \quad 
\Exp_j[\sigma_i(b/2)] &\geq \exp[r_{ij}b/2] \\ \text{ or } \quad 
\Exp_j[\sigma_i(b/2)] &\geq  b^2 / (4\Exp_j[\ell_i^2(1)]).
\end{align*}
Since these both contradict \eqref{contr} for sufficiently large $b$, the proof is complete. 
\end{IEEEproof}

Finally, we prove Theorem~\ref{thm:CUSUM_control} using the definitions and lemmas used thus far in the appendixes. \\

\begin{IEEEproof}[Proof of Theorem~\ref{thm:CUSUM_control}]
We follow the same method as in the  proof of \cite[Theorem 2]{nikiforov2000simple}. 
Without loss of generality, we  assume that~$b$ is large enough so that Lemma~\ref{lem:nonzero_prob} holds. Indeed, for smaller~$b$, the function~$\phi(b)$ can be trivially selected to make the upper bound greater than one. 
To lighten the notation, we suppress dependence on~$b$ and simply write~$\widehat{\sigma}$,~$\sigma$, and~$\sigma_i$ instead of $\widehat{\sigma}(b)$, $ \sigma(b)$, and $ \sigma_i(b)$, respectively.

From  Lemma~\ref{lem:nonzero_prob} 
and the law of total expectation we obtain:
\begin{align}\label{eq1}
\begin{split}
 &\Pro_{\nu, j} ( \sigma_{i} \leq  \sigma_j \, |\; \sigma > \nu) \\
 &= 1- \Pro_{\nu, j} (  \sigma_{i} >  \sigma_j  \, | \, \sigma > \nu)\\
&= 1-  \frac{\Exp_{\nu, j}[\sigma_{i} -  \sigma_j  \wedge \sigma_i \, |  \, \sigma> \nu]}{\Exp_{\nu, j}[\sigma_{i} -  \sigma_j\, |\,  \sigma > \nu, \sigma_{i}>  \sigma_j ]}. 
\end{split}
\end{align}
On the event $$ \lbrace \sigma > \nu, \; \sigma_i  > \sigma_j  \rbrace \in \mathcal{F}_{\sigma_j}$$ we have  $$\nu < \sigma \leq \sigma_j < \sigma_i,$$
and the difference  ${\sigma_i -\sigma_j}$ is a function of $Y_{i}(\sigma_j)$ and 
${ \{X_{n}, n > \sigma_j\} }$.  Recalling the definition of $L_{ij}$ in \eqref{L} we then have 
\begin{align*} 
&\Exp_{\nu, j} \left[ \sigma_i-  \sigma_j \, | \, \sigma > \nu, \sigma_{i} >  \sigma_j, \cF_{\sigma_j}  \right] \\
&= L_{ij}( Y_{i}(\sigma_j) \,;\,b) \; \cdot \; \mathds{1}(\{  \sigma > \nu, \sigma_{i} >  \sigma_j\})\\
& \leq L_{ij}(0;b) = \Exp_j[\sigma_{i}],
\end{align*}
where the inequality holds because~$L_{ij}(\,\cdot \, ; \, b)$ is a non-increasing function (note that  $ Y_i(\sigma_j) < b $ on $ \lbrace \sigma_i > \sigma_j \rbrace$), and the final equality follows by definition. 
 By the law of iterated expectation  and the monotonicity of  expectation  it then  follows that 
\begin{align}\label{ineq2}
 \Exp_{\nu, j} & \left[\sigma_i-\sigma_j \, | \, \sigma > \nu, \sigma_{i} >  \sigma_j \right]  \leq \Exp_j[\sigma_{i}].
\end{align}
Combining~\eqref{eq1}--\eqref{ineq2}, we conclude that: 
 \begin{align*} 
\Pro_{\nu, j} ( \sigma_{i} \leq  \sigma_j \, |\, \sigma> \nu) 
 &\leq 1 -  \frac{\Exp_{\nu, j}[\sigma_i- \sigma_j  \wedge \sigma_i\,|\, \sigma > \nu] }{\Exp_j[\sigma_i]}  ,
 \end{align*}   
and  adding and subtracting~$\nu$ in the numerator of the fraction in the upper bound we further  have 
  \begin{align}  \label{ineq3}
  \begin{split}
\Pro_{\nu, j} ( \sigma_{i}\leq  \sigma_j \, |\, \sigma > \nu) 
\leq   
 \frac{ \Exp_j[\sigma_i] - \Exp_{\nu, j}[\sigma_i - \nu\,|\, \sigma > \nu]}{\Exp_j[\sigma_i ]} \,  \\
 +  \, \frac{ \Exp_{\nu, j}[  \sigma_j \wedge \sigma_i - \nu\,|\, \sigma > \nu]}{\Exp_j[\sigma_{i}] }.
 \end{split}
 \end{align} 
For the first term in the upper bound, we observe that on the event $\{\sigma>\nu\}$ we have  $Y_{i}(\nu) <b$ and that $\sigma_i -\nu$  is  a positive-valued  function of $ Y_i(\nu)$ and $ X_{\nu + 1}, X_{\nu + 2}, \ldots$.  As a result,
 \begin{align}
 \begin{split}
& \Exp_{\nu,j}[\sigma_i-\nu \, | \, \cF_\nu, \sigma>\nu] \\
&=  L_{ij}(Y_{i}(\nu)  \wedge b \, ; \, b) \cdot  \mathds{1} (\{ \sigma>\nu\}) ,
 \end{split}
 \end{align}
 and   by another application of the law of iterated expectation we obtain 
 \begin{align}\label{eqn:total_exp}
\Exp_{\nu, j}[\sigma_{i}- \nu \, | \,  \sigma> \nu] 
= \Exp_{\infty} \left[ L_{ij}(Y_{i}(\nu) \wedge b \,  ; \, b) \, | \, \sigma  >\nu  \right]  ,
\end{align}
where in the equality we have also used the  fact that 
$Y_i(\nu)$ has  the same distribution under $\Pro_{\nu, j}$ and $\Pro_{\infty}$. For the second second term of~\eqref{ineq3}, we observe that 
$$ \Exp_{\nu, j}[  \sigma_j \wedge \sigma_i - \nu \, |\, \sigma > \nu] \leq  \Exp_{\nu, j}[  \sigma_j - \nu \, |\,  \sigma > \nu] \leq \Exp_j[\sigma_j]. $$ 
 Thus, we conclude that 
  \begin{align}  \label{ineq3a}
  \begin{split}
&  \Pro_{\nu, j} ( \sigma_{i}\leq  \sigma_j \, |\, \sigma > \nu) 
   \\
 & \leq   \frac{\Exp_j[ \sigma_i]   - \Exp_{\infty} \left[ L_{ij}(Y_{i}(\nu) \wedge b \, ; \, b) \, | \, \sigma > \nu \right]  }{\Exp_j[\sigma_{i}] }
+ \frac{\Exp_j[ \sigma_j]}{\Exp_j[\sigma_{i}] }.
\end{split}
 \end{align} 
 We bound each of these two terms separately in (i) and (ii). \\


\noindent \textbf{Proof of part (i)}
Suppose ${I_{ji} > I_j}$, and that $\psi_{ij}$ has a positive root.   We first bound the second term of~\eqref{ineq3a}. From~\eqref{eqn:cusum_delay_upper_bound} and~\eqref{eqn:exponential_lower_bound} we obtain
\begin{equation}\label{part2}
     \frac{ \Exp_{j}[  \sigma_j]}{\Exp_j[\sigma_{i}] } \leq e^{-r_{ij}b} \left( \frac{b}{I_j} +  B_j \right),
\end{equation}
thus, it remains to upper bound the  first term in \eqref{ineq3a}. For this, we recall the definitions of $L_{ij}$ and  $l_{ij}$ in \eqref{L}  and~\eqref{def: l}. 
From Lemma~\ref{L_Lemma} we then have 
\begin{align}\label{eqn:split1}
\begin{split}
 &\Exp_j[ \sigma_i] - \Exp_{\infty} \left[ L_{ij}(Y_{i}(\nu) \wedge b \,  ; \, b) \, | \, \sigma  >\nu  \right] \\ 
 &\leq 
   L_{ij}(0\,;\,b) - l_{ij}(0\,;\,b)  \\
   &+  \Exp_{\infty} \left[ l_{ij}(0\,;\,b) -  l_{ij}(Y_{i}(\nu) \wedge b \, ; \, b) \, | \, \sigma > \nu \right].
   \end{split}
\end{align}
From~\eqref{def: l} it follows that
\begin{equation}\label{eqn:first_part1}
L_{ij}(0\,;\,b) - l_{ij}(0\,;\,b) = e^{-r_{ij}b} \left[ \frac{b + \omega_{ij}}{I_{ji} - I_{j}} + L_{ij}(0\,;\,b) \right]. 
\end{equation} 

To bound the second term in the upper bound in~\eqref{eqn:split1}, we proceed as follows:
\begin{align}\label{eqn:FTC}
\begin{split}
&\Exp_{\infty} \left[ l_{ij}(0\,;\,b) -  l_{ij}(Y_{i}(\nu)\wedge b \, ; \, b) \,  | \, \sigma > \nu \right] \\
&= \Exp_{\infty} \left[ - \int_0^b l_{ij}'(x\,;\,b)\,  \mathds{1} (\lbrace Y_i(\nu) \geq x \rbrace) dx \, \Big| \, \sigma > \nu \right] \\
&= - \int_0^b l_{ij}'(x\,;\,b) \,\Pro_{\infty}(Y_i(\nu) \geq x \, | \, \sigma > \nu) \; dx \\
&\leq - \int_0^b l_{ij}'(x\,;\,b) \, e^{-x} \; dx \\
&\leq \left( \frac{b + \omega_{ij}}{I_{ji} - I_j} + L_{ij}(0;b) \right) r_{ij} e^{-r_{ij}b} \int_0^b e^{x(r_{ij} - 1)} dx.
\end{split}
\end{align}
The first equality follows from the  Fundamental Theorem of Calculus, the second by Tonelli's theorem,  the first inequality from  the assumption that  \eqref{condition on change points} holds, and the third one from the bound on $-l'_{ij}$ in Lemma~\ref{lemma:derivatives}. The integral in the upper bound of~\eqref{eqn:FTC}  is upper bounded by 
 \begin{equation}\label{eqn}
\begin{cases} 
 b,  \quad &\text{if} \quad  r_{ij} \leq 1 \\
\frac{1}{r_{ij} - 1} e^{-(1-r_{ij}) b} , \quad & \text{if} \quad  r_{ij} >1.
 \end{cases}
\end{equation} 
In view of this, we distinguish two cases depending on whether   $r_{ij} > 1$ or $r_{ij} \leq 1$. \\

When $r_{ij} >1$,  from 
\eqref{eqn:split1}--\eqref{eqn} it follows that the numerator of the first term in \eqref{ineq3a}  is upper bounded by 
\begin{align}
    \left( \frac{b + \omega_{ij}}{I_{ji} - I_j} + L_{ij}(0\, ;\, b) \right) \cdot \Big( e^{-r_{ij}b} + \frac{r_{ij}}{r_{ij}-1}e^{-b} \Big) 
\end{align}
and, as a result, in view of~\eqref{eqn:exponential_lower_bound}, and the fact that $L_ij(0;b) \equiv \Exp_j[\sigma_i(b)]$, the first term in 
\eqref{ineq3a}  is upper bounded by 
 \begin{align}\label{part1_rlarge}
     \left( e^{-r_{ij}b} \, \frac{b + \omega_{ij}}{I_{ji} - I_j} + 1 \right)  \cdot  \left( e^{-r_{ij}b}  + e^{-b}\frac{r_{ij}}{r_{ij}-1} \right).
 \end{align}


Combining this with \eqref{ineq3a}  and \eqref{part2},  we conclude that
 \begin{equation}
 \Pro_{\nu, j}(\sigma_i \leq \sigma_j \, | \, \sigma > \nu) \leq e^{-b}\frac{r_{ij}}{r_{ij}-1}  (1 + \phi_{ij}(b) ),
    \end{equation}
    where
    \begin{align*}
    \phi_{ij}(b) &\equiv e^{-b(r_{ij}-1)} \frac{r_{ij}-1}{r_{ij}}\, \left(   \frac{b}{I_{j}} + B_j + 1 \right) \\ 
&\quad +    \left( \frac{b + \omega_{ij}}{I_{ji} - I_j} 
 \right) \left( e^{-r_{ij}b} + \frac{r_{ij} - 1}{r_{ij}} e^{-2r_{ij}b + b} \right),
\end{align*}
 which is a positive function that goes to zero as $b \to \infty$ and does not depend on~$\nu$. \\

When $r_{ij} \leq 1$,    from 
\eqref{eqn:split1}--\eqref{eqn} it follows that the numerator of the first term in \eqref{ineq3a}  is upper bounded by
$$    \left( \frac{b + \omega_{ij}}{I_{ji} - I_j} + L_{ij}(0\, ;\, b) \right)\,  e^{-r_{ij}b}\, \Big( r_{ij}b +1 \Big) 
$$
and, as a result, in view of~\eqref{eqn:exponential_lower_bound}, the first term in 
\eqref{ineq3a}  is upper bounded by 
$$
     \left( e^{-r_{ij}b} \, \frac{b + \omega_{ij}}{I_{ji} - I_j} + 1 \right)     e^{-r_{ij}b}\, \Big( r_{ij} b +1 \Big).
$$
Combining this with \eqref{ineq3a}  and \eqref{part2},   we conclude that 
$$
 \Pro_{\nu, j}(\sigma_i \leq \sigma_j \, | \, \sigma > \nu) \leq \\
 b \, e^{-r_{ij}b} \left( r_{ij} + \frac{1}{I_j } \right) \, (1 + \zeta_{ij}(b) ),
$$
    where~$\zeta_{ij}(b)$ is a positive function that goes to zero as $b \to \infty$ which does not depend on~$\nu.$\\

\textbf{Proof of part (ii)}
We now assume that $\Exp_j[\ell_i(1)] = 0$\@ and that $\psi_{ij}$ is finite around 0. As before,  we bound each term in the upper bound  of~\eqref{ineq3a}. For the second,  from~\eqref{eqn:quadratic_lower_bound}  and  \eqref{eqn:cusum_delay_upper_bound} we obtain
\begin{equation}
\begin{split}\label{eqn:last_part4}
\frac{\Exp_{ j}[ \sigma_j]}{\Exp_j[\sigma_i]} 
&\leq \frac{\Exp_j[\ell_i^2(1)]}{b^2} \, \left( \frac{b}{I_j}  +  B_j \right).
\end{split}
\end{equation}
For the first one, we recall the definition of $ u_{ij}$ from~\eqref{def: u}. Then, from Lemma~\ref{L_Lemma} it follows that  
\begin{align}\label{eqn:split2}
\begin{split}
& \Exp_j[ \sigma_i] - \Exp_{\infty} \left[ L_{ij}(Y_{i}(\nu) \wedge b \,  ; \, b) \, | \, \sigma  >\nu  \right] \\ 
 &\leq     L_{ij}(0\,;\,b) - u_{ij}(0\,;\,b)   \\
   &\quad +  \Exp_{\infty} \left[ u_{ij}(0\,;\,b) -  u_{ij}(Y_{i}(\nu) \wedge b \, ; \, b) \, | \, \sigma > \nu \right] .
  \end{split}
\end{align}
By~\eqref{def: u} we have:
\begin{equation}\label{eqn:diff2}
   L_{ij}(0\,;\,b) - u_{ij}(0\,;\,b)
   = L_{ij}(0\,;\,b) \, \frac{\omega_{ij} + \tilde{\omega}_{ij}}{b + \omega_{ij} + \tilde{\omega}_{ij}}.
\end{equation}
Moreover, we can  upper bound the second term in~\eqref{eqn:split2} as follows:
\begin{align}\label{old_3}
\begin{split}
& \Exp_{\infty} \left[ u_{ij}(0\,;\,b) -  u_{ij}(Y_{i}(\nu) \wedge b \,;\, b) \,|\, \sigma > \nu \right] \\
& =  \Exp_{\infty} \left[ \int_0^b -u'_{ij}(x\,;\,b) \mathds{1} ( \lbrace Y_i(\nu) \geq x \rbrace )  dx \, \big| \, \sigma > \nu \right]  \\
 & = \int_0^b -u'_{ij}(x\,;\,b) \,  \Pro_{\infty} ( Y_i(\nu) \geq x \,\big| \,\sigma > \nu) dx \\
  & \leq  \int_0^b -u'_{ij}(x\,;\,b) \,   e^{-x} \; dx \\
  &\leq \frac{L_{ij}(0;b)}{b + \omega_{ij} + \tilde{\omega}_{ij}}  +\frac{b}{\Exp_j[\ell_i(1)^2] \, (b + \omega_{ij} + \tilde{\omega}_{ij})}.
 \end{split}
\end{align}
The first equality follows by the  Fundamental Theorem of Calculus and the second by Tonelli's theorem. The first inequality follows by the assumption that \eqref{condition on change points} holds, and the second  by the upper bound on $-u'_{ij}$ in Lemma~\ref{lemma:derivatives}\@. From~\eqref{eqn:split2}--\eqref{old_3} and~\eqref{eqn:quadratic_lower_bound} , and the fact that $L_{ij}(0;b) \equiv \Exp_j[\sigma_i]$, it then follows that the
first term in \eqref{ineq3a} is upper bounded by 
\begin{align*}\label{eqn:last_part5}
& \frac{\omega_{ij} + \tilde{\omega}_{ij}}
{ b + \omega_{ij} + \tilde{\omega}_{ij}} + \frac{1}{b(b + \omega_{ij} + \tilde{\omega}_{ij})} \\
&\leq  \frac{(1 + \omega_{ij} + \tilde{\omega}_{ij})}{b} + \frac{1}{b^2} ,
\end{align*}
where the last inequality follows from the fact that~$\omega_{ij}$ and~$\tilde{\omega}_{ij}$ are non-negative.  Combining this with  \eqref{ineq3a}  and \eqref{eqn:last_part4}, we arrive at
\begin{align*}
    \Pro_{\nu, j}(\sigma_i \leq \sigma_j \, | \, \sigma > \nu) & \leq
\widetilde{C} \, b^{-1}  \, \left( 1 + \xi_{ij}(b) \right)  ,
\end{align*}
where 
\begin{align*}
\widetilde{C} &\equiv  1 + \omega_{ij} + \tilde{\omega}_{ij} + \frac{\Exp_j[\ell_i(1)^2]}{I_j}  \\
    \xi_{ij}(b) & \equiv
    \left(B_j \, \Exp_j[\ell_i^2(1)] + 1 \right)/ (\widetilde{C} \, b).
\end{align*}
The latter is a positive function that goes to zero as~$b \to \infty$ and does not depend on~$\nu$, so this completes the proof. 
\end{IEEEproof}

\section{Stochastic Monotonicity of the CuSum}\label{app:lemmas}

\begin{lemma}\label{lemma:Multichannel_assumption}
For all $i \in \mathcal{I}$, the CuSum process ${ \lbrace Y_i (n), n \in \mathbb{N} \rbrace}$ is a non-negative stochastically monotone Markov process under~$\Pro_\infty .$
\end{lemma}
\begin{IEEEproof}
Specifically, we show that $${\Pro_{\infty}(Y_i(1) \geq y \, | \, Y_i(0) = x)}$$ is non-decreasing and right-continuous in~$x$ for all~$y \geq 0$\@.

By the definition of the CuSum statistic, \begin{equation} \Pro_{\infty}(Y_i(1) \geq y \, | \, Y_i(0) = x) = \Pro_{\infty}( (\ell_i(1) + x)^+ \geq y).
\end{equation} For $y=0$, \begin{equation}
\Pro_\infty( (\ell_i(1) + x)^+ \geq y) =1 \, \quad \text{ for all } x.
\end{equation}
On the other hand, for $y>0,$ 
\begin{equation}
 \Pro_\infty( (\ell_i(1) + x)^+ \geq y)  
    = \Pro_\infty( \ell_i(1) \geq y -x)  ,
\end{equation}
which is non-decreasing in~$x$ for all~$y > 0$\@.


Let $ \epsilon > 0$ and fix $ x \geq  0$ and $y > 0$\@. (When $y=0$, the probability is constant with respect to $x$.) We show right-continuity:
\begin{align*}
&\Pro_{\infty}(Y_i(1) \geq y \, | \, Y_i(0) = x + \epsilon) - \Pro_{\infty}(Y_i(1) \geq y \, | \, Y_i(0) = x ) \\
&= \Pro_{\infty}( \ell_i(1) \geq y - (x + \epsilon)) - \Pro_{\infty}( \ell_i(1) \geq y - x)  \\
&= \Pro_{\infty} ( (y-x) - \epsilon \leq \ell_i(1) < (y-x)),
\end{align*} 
but $$ \Pro_{\infty} ( (y-x) - \epsilon \leq \ell_i(1) < (y-x)) \to 0 \text{ as } \epsilon \to 0. $$
\end{IEEEproof}


\ifCLASSOPTIONcaptionsoff
  \newpage
\fi

\end{document}